\documentclass[12pt]{article}
\usepackage{latexsym}
\newtheorem{thm}{Theorem}
\newtheorem{con}{Conjecture}
\begin{document}
\title{\bf Complete minors and stability numbers\thanks{Project 11771246 supported by National Natural Science Foundation of China.}}
\author{{\sc Wenkai Fu and Lingsheng Shi}\\ Department of Mathematical Sciences\\ Tsinghua
University, Beijing, 100084, China\\ E-mail:
lshi@math.tsinghua.edu.cn}
\date{}
\maketitle
\begin{abstract}
Hadwiger's conjecture implies that $n\le\alpha h$ for all graphs of order $n$, stability number $\alpha $, and Hadwiger number $h$. Combining ideas of Kawarabayashi et al. and Wood, we prove that $n\le (\alpha-1)(2h-5)+5$ for such graphs if $\alpha\ge 3$ and $h\ge 5$.

\end{abstract}

Let $G=(V,E)$ be a finite simple graph. The {\em order} of $G$ is $|V|$, the number of vertices. A {\em stable set} of $G$ is a subset of vertices which are pairwise nonadjacent. The {\em stability number} of $G$, denoted by $\alpha (G)$, is the size of the largest stable set. A {\em clique} of $G$ is a complete subgraph of $G$. The {\em clique number} of $G$, denoted by $\omega (G)$, is the order of the largest clique in $G$. The {\em chromatic number} of $G$, denoted by $\chi (G)$, is the smallest integer $k$ such that the vertex set $V$ can be partitioned into $k$ stable sets. A complete graph is called a {\em complete minor} of $G$ if it can be obtained by contracting a subgraph of $G$. The {\em Hadwiger number} of $G$, denoted by $h(G)$, is the order of its largest complete minor.

In 1943, Hadwiger \cite{H} came up with a conjecture which generalized the four color theorem \cite{AH,AHK}. It has been widely considered as one of the most interesting and important problems in graph theory, see \cite{S} for a survey.

\medskip\noindent
{\bf Hadwiger's Conjecture}. {\em For every graph $G$, $\chi (G)\le h(G)$}.

\medskip
Hadwiger \cite{H} proved his own conjecture for $\chi\le 4$. Wagner \cite{W} and Robertson et al. \cite{RST} proved the equivalence of the Hadwiger conjecture and the four color theorem for $\chi =5$ and $\chi =6$, respectively. This conjecture is still open for $\chi\ge 7$. Since it is obvious that $|V(G)|\le\alpha (G)\chi (G)$, the Hadwiger conjecture implies the following result.
\begin{con}\label{wd}
If $G$ is a graph of order $n$ with stability number $\alpha $ and Hadwiger number $h$, then $n\le\alpha h$.
\end{con}
Conjecture~\ref{wd} seems weaker than the Hadwiger conjecture, however for $\alpha =2$ the two conjectures are equivalent which is proved by Plummer et al. \cite{PST}. Conjecture~\ref{wd} holds for $h\le 5$ since the Hadwiger conjecture holds for $\chi\le 6$. Though Conjecture~\ref{wd} was explicitly stated by Woodall \cite{Wo} in 1987, it had been studied before its publication. In fact, the first weak version of Conjecture~\ref{wd} was obtained in 1982 by Duchet and Meyniel \cite{DM} who proved that
$$n\le (2\alpha -1)h.$$
There have been several improvements on their result. In 2005, Kawarabayashi et al. \cite{KPT} proved that
\begin{equation}\label{kp}
n\le(2\alpha -1)(h-1)+1
\end{equation}
and $n\le (2\alpha -1)h-\omega $ for $\alpha\ge 2$ and $n\le (2\alpha -3/2)h$ for $\alpha\ge 3$, which was further improved by Kawarabayashi and Song \cite{KS} to
\begin{equation}\label{ka}
n\le 2(\alpha -1)h.
\end{equation}
In 2007, Wood \cite{Wood} came up with another improvement on (\ref{kp}) by showing that
\begin{equation}\label{wh}
n\le (2\alpha -1)(h-5/2)+5/2\mbox{ for }h\ge 5.
\end{equation}
In 2010, Fox \cite{F} was the first to improve on the factor 2 by proving that $n\le 1.983\alpha h$, which was slightly improved by Balogh and Kostochka \cite{BK} in 2011 to $n\le 1.948\alpha h$.

Combining the ideas of Kawarabayashi et al. \cite{KPT,KS} and Wood \cite{Wood}, we make an improvement on both (\ref{ka}) and (\ref{wh}), which is also better than the bounds of Fox and of Balogh and Kostochka when $\alpha $ or $h$ is small.
\begin{thm}
If $G$ is a graph of order $n$ with stability number $\alpha\ge 3$ and Hadwiger number $h\ge 5$, then $n\le (\alpha -1)(2h-5)+5$.
\end{thm}
{\bf Proof}. We use induction on the Hadwiger number $h$ of $G$. This theorem holds for $h=5$ since Conjecture~\ref{wd} holds for $h\le 5$. Now consider such a graph $G$ of order $n$ with $h=h(G)\ge 6$ and $\alpha =\alpha (G)\ge 3$.

\medskip
{\bf Case 1}. The graph $G$ is disconnected.

\medskip
In this case assume that $G$ is a disjoint union of two nonempty subgraphs, say $G_1$ and $G_2$. Let $n_i=|V(G_i)|$, $\alpha _i=\alpha (G_i)$, and $h_i=h(G_i)$ for $i=1,2$. It is clear that $n=n_1+n_2$, $\alpha =\alpha _1+\alpha _2$, and $h=\max\{ h_1,h_2\}$. We observe that $n_i\le (2\alpha _i-1)(h-5/2)+5/2$ for $i=1,2$. This observation follows from (\ref{wh}) for $h\ge h_i\ge 5$. For $h_i<5<h$, Conjecture~\ref{wd} holds for $G_i$ and $n_i\le\alpha _ih_i<5\alpha _i\le (2\alpha _i-1)(h-5/2)+5/2$. Thus
$$n\le (2\alpha _1-1)(h-5/2)+5/2+(2\alpha _2-1)(h-5/2)+5/2=(\alpha -1)(2h-5)+5.$$

{\bf Case 2.} The graph $G$ is connected.

\medskip
A claw of $G$ is an induced subgraph $K_{1,3}$. Chudnovsky and Frakdin \cite{CF} proved Conjecture~\ref{wd} for claw-free connected graphs with $\alpha\ge 3$. Therefore, if $G$ is claw-free, then $n\le\alpha h<\alpha h+(\alpha -2)(h-5)=(\alpha -1)(2h-5)+5$. Thus $G$ has a claw and then we construct a connected dominating set. Start with a claw $C$ of $G$ and let $D_0=V(C)$. It is obvious that $D_0$ is connected. If $D_i$ does not dominate the whole graph $G$, then there is a vertex $v\in V(G)\setminus D_i$ such that the distance from $v$ to $D_i$ is 2. Let $P$ be such a path of length 2 linking $v$ and $D_i$, and let $u$ be the centre of $P$ which is adjacent to both $v$ and $D_i$. We put the pair of vertices $u$ and $v$ to $D_i$ and obtain $D_{i+1}$. Apparently $D_{i+1}$ is still connected. Repeat the procedure until $D_k$ is dominating for some $k$.

Let $D$ denote our final connected dominating set and let $S$ be a maximum stable set of the subgraph $G[D]$ induced by $D$. Since $k$ pairs of vertices are in total put into $D_0$, we have $|D|=2k+4$. Another observation is the size of the stable set increased by 1 every time a pair of vertices were put in, because the corresponding vertex $v$ is not adjacent to $D_i$. By $\alpha (C)=3$, we have $k+3\le |S|\le\alpha $ and thus $|D|\le 2|S|-2\le 2\alpha -2$.

Let $H=G-D$, $n_0=|V(H)|$, $\alpha _0=\alpha(H)$, and $h_0=h(H)$. Since $D$ is a connected dominating set, we get $h_0\le h-1$.

\medskip
{\bf Claim 1}. $h_0\ge 5$.

\medskip
If $h_0<5$, then as above, $n_0\le\alpha _0h_0<5\alpha _0\le 5\alpha $ and thus
$$n=|D|+n_0<2\alpha -2+5\alpha =7\alpha-2\le (\alpha -1)(2h-5)+5\mbox{ for }h\ge 6.$$

{\bf Claim 2}. $\alpha _0\ge 3$.

\medskip
Suppose $\alpha _0\le 2$. By (\ref{kp}), we have $n_0<3h_0\le 3h-3$ and
$$n=|D|+n_0<2\alpha -2+3h-3=2\alpha +3h-5.$$
It follows that
\begin{eqnarray*}
n-(\alpha -1)(2h-5)-5&< & 2\alpha +3h-5-2\alpha h+5\alpha +2h-10\\
&=& 7\alpha +5h-2\alpha h-15\\
&=& 5/2-(\alpha -5/2)(2h-7)\le 0
\end{eqnarray*}
for $\alpha\ge 3$ and $h\ge 6$. This proves Claim 2.

\medskip
By Claims 1 and 2 and applying the inductive hypothesis to $H$, we get
$$n=|D|+n_0\le 2\alpha -2+(\alpha _0-1)(2h_0-5)+5\le (\alpha -1)(2h-5)+5$$
for $\alpha _0\le\alpha $ and $h_0\le h-1$. \hfill $\Box $

\end{document}